\documentclass{article}
\usepackage{amsthm,amsfonts,amsbsy, amssymb,amsmath,graphicx}

\author{Vassily Olegovich Manturov}

\date{}

\title{Minimal diagrams of classical knots}

\newtheorem{thm}{Theorem}
\newtheorem{lm}{Lemma}
\newtheorem{st}{Statement}
\newtheorem{re}{Remark}

\begin{document}

\maketitle

\abstract{We show that if a classical knot diagram satisfies a
certain combinatorial condition then it is minimal with respect to
the number of classical crossings. This statement is proved by
using the Kauffman bracket and the relation between atoms and
knots.}

\section{The main result}

In paper \cite{MaArx}, we showed that if a virtual link diagram
satisfies two certain conditions (one of them deals with the
Kauffman bracket and the other one uses the Khovanov homology)
then this diagram is minimal with respect to the number of
classical crossings. That result is a generalisation of the famous
Kauffman-Murasugi theorem \cite{Mur}.

In the present paper, we show that the Kauffman bracket itself
(without Khovanov's categorification) is indeed a very strong tool
to establish minimality of knot diagrams. The condition described
in the present paper deals only with some combinatorics of the
knot diagram, namely, with so-called {\em atom}.

This condition is very easy to check unlike that in \cite{MaArx}
where one should be able to calculate (a part of) the Khovanov
homology. The condition of the present paper uses only some simple
combinatorics.

Note that though the techniques of the present paper uses much
from virtual knot theory, the main theorem is stated {\bf only for
classical knots}, i.e., not virtual knots or links and not
classical links. The reason is that one important step of the
proof deals with the connected summation which is not well defined
either for links or for virtual knots.

What remains in the general case of virtual links, is the
analogous {\em framed} result, which was first proved in
\cite{Ma}. As for virtual {\em knots} (not links), one can obtain
a similar result in the {\em long category}. For long virtual
knots, see \cite{MaLong}.

We shall give the first definition for the general virtual case,
however, to understand the main line of the present paper and the
proof of the main theorem one need not know virtual knot theory.

Virtual knots were proposed by Kauffman in \cite{KaV}. All
necessary detailed definitions can be found therein.

The theory of atoms and knots is represented in \cite{Ma}. Recall
the main definitions.

By an {\em atom} we mean a pair $(M,\Gamma)$, where $M$ is a
connected $2$-manifold and $\Gamma$ is a $4$-valent graph together
with a rule for embedding in $M$ such that the complement
$M\backslash \Gamma$ admits a checkerboard colouring. The graph
$\Gamma$ is called {\em the frame} of the atom. We also think that
for a given atom the colouring is fixed.

The atoms are considered up to natural equivalence, i.e.
homeomorphisms mapping the frame to the frame and preserving the
colour of edges. Certainly, the atom (up to equivalence) is
nothing but its frame together with the rule for attaching black
cells at each vertex (the way for attaching white cells is defined
automatically together with the structure of opposite edges at
vertices).

Let $L$ be a virtual link diagram. Let us construct the atom
$V(L)$ as follows. First, we construct the frame $\Gamma$ of
$V(L)$. The vertices of the frame are in one-to-one correspondence
with {\em classical} crossings of the diagram $L$. Classical
crossings are connected by {\em arcs} which might intersect or
have selfintersection at some virtual crossings. In the classical
case there are no other crossings, so the branches are just the
edges of the shadow of the knot (link). Each classical crossing
has four emanating branches. We associate four edges of the atom
to these branches.

Then the rule for attaching black $2$-cells to the frame is
recovered from the diagram $L$. Namely let $X$ be a classical
crossing of $L$. Enumerate the four emanating branches by letters
$x_{1},x_{2},x_{3},x_{4}$ in the clockwise direction in such a way
that the edges $x_{1}$ and $x_{3}$ form an undercrossing, whence
the edges  $x_{2},x_{4}$ form an overcrossing. Then, the {\em
black} angles are chosen to be  $(x_{1},x_{2})$ and
$(x_{3},x_{4})$

Let $L$ be a virtual diagram and let  $V(L)$ be the corresponding
atom. Each vertex of the atom $V(L)$ is incident to four pieces of
cells: two black ones and two white ones. Globally, some of them
(e.g. two white ones) may coincide. A diagram $L$ is said to be
{\em good} if at each vertex of the atom $V(L)$ we have precisely
four different cells.

To be concise, we shall say {\em genus} (or {\em Euler
characteristic}) of the diagram $L$ for the genus (resp., Euler
characteristic) of the corresponding atom:
$\chi(L)=\chi(V(L)),g(L)=g(V(L))$.

The main result of the present paper is the following
\begin{thm}
Suppose the diagram  $L$ of a classical knot is good. Then it is
minimal. In other words, if the diagram $L$ has $n$ crossings then
for any classical diagram representing the same knot the number of
crossing is at least $n$.
\end{thm}

\begin{re}
There is an important conjecture whether a minimal classical
diagram is minimal in the virtual category, i.e. whether there
exists a minimal classical link $L$ with $n$ crossings admitting a
virtual diagram with strictly smaller number of classical
crossings.

 The main theorem of this paper does not give even a
partial answer to this theorem: we say that a classical diagram is
minimal only among classical ones
\end{re}

To prove this theorem, we shall use some auxiliary lemmas.

By $span X$ for a one-variable (Laurent) polynomial $X$ we mean
the difference between its leading degree and lowest degree. For a
polynomial in many variables, we also may define $span$ with
respect to any of these variables.

\begin{lm}
Let $L$ be a virtual link diagram. Then the following inequality
holds:

\begin{equation}
span\langle L\rangle \le 4n+2(\chi(L)-2),\label{nera}
\end{equation}
whence if $L$ is a good diagram then this inequality (\ref{nera})
becomes a strict equality.
\end{lm}

The proof can be found in e.g. \cite{Ma} or \cite{MaArx}.

We shall use the operation of taking $k$ parallel copies of the
link $L\to D_{k}(L)$. This operation is well defined only in the
{\em framed} category. Framed (virtual) links are equivalence
classes of virtual link diagrams modulo (generalised Reidemeister
moves), where we do not allow the first classical Reidemeister
move and replace it by the double twist move (addition/removal of
two loops having opposite signs), for more detais, see,
e.g.,\cite{Ma}.

More precisely, the following lemma holds.

\begin{st}
If  $L,L'$ are equivalent virtual link diagrams so that the writhe
numbers (framings) of the corresponding components for $L$ and
$L'$ coincide then for each natural  $m$, the diagrams
$D_{m}(L),D_{m}(L')$ represent equivalent virtual links.
\end{st}

\begin{lm}
Suppose a virtual link diagram $L$ is good. Then for any natural
$k$, the diagram $D_{k}(L)$ is good as well.
\end{lm}

The proof can be found in \cite{Ma}. The main idea is that to each
cell of the atom $V(L)$, there correspond precisely $k$
``parallel'' cells of the atom $V(D_{k}(L))$. If a cell
$C_{1},C_{2}$ of $V(D_{k}(L))$ touches itself at a crossing $X$ of
$D_{k}(L)$ then the the corresponding cell in $L$ touches itself
at the crossing corresponding to $X$: for each crossing of $L$, we
have $k^{2}$ corresponding crossings of $D_{k}(L)$.

For any virtual link diagram $L$, its {\em mirror} diagram ${\bar
L}$ is defined to be the diagram obtained from $L$ by switching
all classical crossings (overcrossings are replaced by
undercrossings and vice versa).

Obviously, the atom $V({\bar L})$ is obtained from $V(L)$ by
changing the colour of the cells. Thus, if $L$ is a good diagram,
then so is ${\bar L}$.

We shall also use the notion of  {\em connected sum} $K_{1}\#
K_{2}$ for two oriented classical knots or two oriented long
virtual knots. Note that the connected sum is not well defined for
links; it is not well defined for {\em compact} virtual knots,
either: it depends on the choice of break points.

Nevertheless, for any two virtual link diagrams $L_{1}$ and
$L_{2}$ we can take {\em any} of its connected sums. We shall use
the notation $K_{1}{\#}K_{2}$ only for the classical connected
sum, which is well defined. In this case, the following lemma
holds.

\begin{lm}
Let  $K$ be a good diagram of a virtual link. Then the diagram
$K{\#}{\bar K}$ is good as well.
\end{lm}

\begin{proof}
Assume the contrary. Suppose that  $l$ is a cell of the atom
$V(K{\#}{\bar K})$ (say, black) that touches itself at some
crossing $X$. The boundary $\partial l$ of this cell is a cycle on
the frame of the atom. On the atom $V(K{\#}{\bar K})$, we have two
edges $e_{1},e_{2}$, ``separating'' $V(K_{1})$ from $V(K_{2})$.
Choose points $X_{1},X_{2}$ on these edges. These points divide
$\partial l$ into two parts. One part of the cycle $\partial l$
generates a black cell of the atom $V(K)$ whence the other one
generates a black cell for $V({\bar K})$. By definition, the
diagram containing the vertex $X$ ($K$ or ${\bar K}$) is not good.
This means that the diagram $K$ is not good. The contradiction
completes the proof.
\end{proof}

The following fact is evident
\begin{st}
If two classical knot diagrams $K_{1},K_{2}$ are isotopic, then
the diagrams $K_{1}\# {\bar K}_{1}$ and $K_{2}\#{\bar K}_{2}$ are
framed equivalent.\label{frq}
\end{st}

\begin{proof}
Indeed, the two knots in question are isotopic and have the same
framing equal to zero.
\end{proof}

Let us now prove the main theorem of this paper. Let  $K$ be a
good classical knot diagram having $n$ crossings. Denote the Euler
characteristic of the diagram $K{\#}{\bar K}$ by $\chi$. Set
$N=2n$. Suppose a classical diagram $K'$ having $n'$ crossings
($n'<n$) generates the same knot as $K$. Denote the genus of the
diagram $K'{\#}{\bar K'}$ by $\chi'$.

Let $m$ be a positive integer. It follows from Statement \ref{frq}
that $D_{m}(K{\#}{\bar K})$ and $D_{m}(K'{\#}{\bar K}')$ generate
isotopic knots. Denote $D_{m}(K{\#}{\bar K})$ by $D_{m}$ and
denote $D_{m}(K'{\#}{\bar K}')$ by $D'_{m}$. By definition we have
$\langle D_{m}\rangle=\langle D'_{m}\rangle$. Also, set
$\chi=\chi(K{\#}{\bar K}),\chi'=\chi(K'{\#}{\bar K'})$.

The diagram $K{\#}{\bar K}$ is good. Thus, so are all diagrams
$D_{m}$ for arbitrary positive integers $m$:

\begin{equation}
span \langle D_{m}\rangle =4 m^{2}N+2(\chi_{m}-2),
\end{equation}

where $\chi_{m}=\chi(D_{m})$. The atom $V(D_{m})$ has $m^{2}N$
vertices,  $2m^{2}N$ edges and $m\Gamma$ $2$-cells, where
$\Gamma=N+\chi$ is the number of the $2$-cells of the atom
$K{\#}{\bar K}$. Thus,

\begin{equation}
span \langle D_{m}\rangle =2 (m^2+m)N+2m \chi-4.\label{ra1}
\end{equation}

Analogously, for the diagram  $D'_{m}$ the following inequality
holds

\begin{equation}
span \langle D'_{m}\rangle \le 2 (m^2+m)N'+2m \chi'-4.\label{ra2}
\end{equation}

Here we can not say whether the exact equality takes place, since
we do not know whether the diagram $D'_{m}$ is good.

Comparing the right-hand sides of (\ref{ra1}) and (\ref{ra2}), we
get

\begin{equation}
(\chi'-\chi)\ge (m+1)(N-N').
\end{equation}

According to the assumption, we have $n-n'>0$; thus $N-N'>0$.
Since  $m$ is chosen arbitrarily, we get to a contradiction: the
quantity $\chi'-\chi$ (which is fixed and does not depend on $m$)
should exceed any preassigned positive integer number. The
contradiction completes the proof of the theorem.

\section{The general case of virtual links}

The proof given above works neither for classical links nor for
virtual knots because of the connected summation. The trick using
the connected sum is indeed needed to compare the diagrams
$D_{m}(K_{1})$ and $D_{m}(K_{2})$. But since we do not know
whether framings of the knot $K_{1}$ (whose minimality is being
tested) and $K_{2}$ coincide, we can not say whether
$D_{m}(K_{1})$ and $D_{m}(K_{2})$ generate isotopic links. In
order to avoid the problem with framing, we have to take the
connected sum of the initial knot with its inverse, thus
restricting ourselves only for the case of classical knots.

However, the trick using the connected sum with the inverse image
is unnecessary, if we deal with {\em framed} knots. This leads to
the following

\begin{thm}
Let $L_{1}$ be a good diagram of a framed virtual link. Then it is
minimal in the framed category.
\end{thm}

This theorem was proved in \cite{Ma}.

Also, if we deal with {\em long} virtual knots (i.e. virtual knots
with fixed endpoints), we have a well-defined connected sum
operation.

This leads to the following

\begin{thm}
Let $K_{1}$ be a long virtual knot diagram such that the
corresponding compact virtual knot diagram $Cl(K_{1})$ is good.
Then the diagram $K_{1}$ is minimal in the long category.
\end{thm}

The proof literally repeats the proof of the main theorem of this
paper.

\section{Examples}

Actually, it is not difficult to construct a non-alternating
classical knot whose minimality can be detected by the main
theorem of this paper. For instance, so are knots represented by
closures of positive braids with arbitrary number of strands,
where we use only exponents of standard generators
$\sigma_{i}^{j}$ for $j\ge 2$ (one should just make sure that the
obtained diagram is a knot, not link).

\begin{re}
Note that if we can apply a Reidemeister move decreasing the
number of crossings to a classical knot diagram $L$ then the
condition of the theorem evidently fails. Thus, we do not lose
generality: evidently non-minimal diagrams could not be minimal.
However, if we can apply a third Reidemeister move to a diagram
$L$, then the diagram does not satisfy the condition of the main
theorem. Thus, this theorem works only for ``fixed'' diagrams,
i.e. those for which any Reidemeister move should increase the
number of crossings.
\end{re}

\end{document}